\documentclass[10pt,a4paper,draft]{amsart}
\usepackage{amssymb}
\def\CC{\mathbb C}

\def\NN{\mathbb N}

\def\RR{\mathbb R}
\def\ZZ{\mathbb Z}

\def\CO{\mathcal O}

\def\bs{\boldsymbol}
\def\Cal{\mathcal}

\let\epsilon=\varepsilon
\def\eps{\epsilon}
\def\Dqed{\text{\popQED}\tag*{\qed}}
\def\halfskip{\vskip 10pt plus 1pt minus 1pt}

\def\ord{\operatorname{ord}}
\let\phi=\varphi
\def\PSH{\mathcal P\mathcal S\mathcal H}

\def\tfrac#1#2{{\textstyle\frac#1#2}}
\def\theequation{\arabic{chapter}.\arabic{section}.\arabic{equation}}
\def\too{\longrightarrow}
\def\wdht{\widehat}
\def\wdtl{\widetilde}
\def\Delta{\varDelta}
\def\Gamma{\varGamma}
\def\Omega{\varOmega}
\def\Phi{\varPhi}
\def\Psi{\varPsi}

\catcode`@=11
\def\@opargbegintheorem#1#2#3{%
\item[\hskip\labelsep \theorem@headerfont #1\ #2]\rom{(#3)}.}
\def\th@mytheorem{%
  \let\thm@indent\noindent 
  \thm@headfont{\bfseries}
  \itshape 
}
\def\th@myremark{%
  \let\thm@indent\noindent 
  \thm@headfont{\bfseries}
}
\catcode`@=12

\theoremstyle{mytheorem}
\newtheorem{Theorem}{Theorem}[section]
\newtheorem{Corollary}[Theorem]{Corollary}
\newtheorem{Lemma}[Theorem]{Lemma}
\newtheorem{Proposition}[Theorem]{Proposition}
\newtheorem{ETh}[Theorem]{}
\theoremstyle{myremark}
\newtheorem{Definition}[Theorem]{Definition}

\newtheorem{Example}[Theorem]{Example}
\catcode`@=11
\tagsleft@false
\def\@eqnnum{\tagform@\theequation}
\def\@makefnmark{\hbox{$\left(^{\@thefnmark}\right)$\;}}
\@mparswitchfalse
\renewcommand\qed{\RIfM@\else\unskip\nobreak\fi\hfill\qedsymbol}
\def\partrunhead#1#2#3{\small%
  \@ifnotempty{#2}{{\ignorespaces#1 #2\unskip}\@ifnotempty{#3}{. }}%
  \def\@tempa{#3}%
  \ifx\@empty\@tempa\else\@tempa\fi}

\def\@cite#1#2{{%
 \m@th\upshape\mdseries[{#1}{\if@tempswa, #2\fi}]}}
\@ifundefined{cite }{%
  \expandafter\let\csname cite \endcsname\cite
  \edef\cite{\@nx\protect\@xp\@nx\csname cite \endcsname}%
}{}
\catcode`@=12

\begin{document}

\title{On extremal holomorphically contractible families}

\author{Marek Jarnicki}
\address{Jagiellonian University, Institute of Mathematics\newline
\indent Reymonta 4, 30-059 Krak\'ow, Poland}
\email{jarnicki@im.uj.edu.pl}
\thanks{The first and the third author were supported in part by KBN grant
no.~5 P03A 033 21 and the Nieders\"achsisches
Ministerium f\"ur Wissenschaft und Kultur, Az. 15.3 -- 50 113(55) PL}

\author{Witold Jarnicki}
\address{Jagiellonian University, Institute of Mathematics\newline
\indent Reymonta 4, 30-059 Krak\'ow, Poland}
\email{wmj@im.uj.edu.pl}
\thanks{The second author was supported in part by KBN grant
no.~2 P03A 015 22.}

\author{Peter Pflug}
\address{Carl von Ossietzky Universit\"at Oldenburg, Fachbereich Mathematik\newline
\indent Postfach 2503, D-26111 Oldenburg, Germany}
\email{pflug@{mathematik.uni-oldenburg.de}}

\subjclass{32F45, 32U35}


\begin{abstract}
We prove (Theorem \ref{T1}) that the category of generalized holomorphically contractible
families (Definition \ref{Def}) possesses maximal and minimal objects. Moreover,
we present basic properties of these extremal families.
\end{abstract}

\maketitle


\section{Introduction. Main results.}

First recall the standard definition of a holomorphically contractible family
(cf.~\cite{JarPfl1993}, \S\;4.1). A family $(d_G)_G$ of
functions
$$
d_G:G\times G\too\RR_+:=[0,+\infty),
$$
where $G$ runs over all
domains $G\subset\CC^n$ with arbitrary $n\in\NN$, is said to be {\it
holomorphically contractible} if the following two conditions are satisfied:

$\bullet$ for the unit disc $E$ we get
$d_E(a,z)=m_E(a,z):=|\frac{z-a}{1-\overline az}|$, $a,z\in E$,

$\bullet$ for any domains $G\subset\CC^n$, $D\subset\CC^m$, every holomorphic
mapping $F:G\too D$ is a contraction with respect to $d_G$ and $d_D$, i.e.
$d_D(F(a),F(z))\leq d_G(a,z)$, $a,z\in G$.

\halfskip

Let us recall some important holomorphically contractible families:

$\bullet$ {\it M\"obius pseudodistance:}
\begin{align*}
c_G^\ast(a,z):=&\sup\{m_E(f(a),f(z)): f\in\CO(G,E)\}\\
=&\sup\{|f(z)|: f\in\CO(G,E),\;f(a)=0\},
\end{align*}

$\bullet$ {\it higher order M\"obius function:}
$$
m_G^{(k)}(a,z):=\sup\{|f(z)|^{1/k}: f\in\CO(G,E),\;\ord_a\!f\geq k\},\quad
k\in\NN,
$$
where $\ord_a\!f$ denotes the order of zero of $f$ at $a$,

$\bullet$ {\it pluricomplex Green function:}
\begin{multline*}
g_G(a,z):=\sup\{u(z):\; \; u:G\too[0,1),\; \log u\in\PSH(G),\\
\exists_{C=C(u)>0}\;\forall_{w\in G}:\; u(w)\leq C\|w-a\|\},
\end{multline*}
where $\PSH(G)$ denotes the family of all functions plurisubharmonic on $G$,

$\bullet$ {\it Lempert function:}
\begin{align*}
\wdtl k_G^\ast(a,z):=&\inf\{m_E(\lambda,\mu): \exists_{\phi\in\CO(E,G)}:\;
\phi(\lambda)=a,\;\phi(\mu)=z\}\\
=&\inf\{|\mu|: \exists_{\phi\in\CO(E,G)}:\;\phi(0)=a,\;\phi(\mu)=z\}.
\end{align*}

It is well known that
$$
c_G^\ast=m_G^{(1)}\leq m_G^{(k)}\leq g_G\leq\wdtl k_G^\ast,
$$
and for any holomorphically contractible family $(d_G)_G$ we have
\begin{gather*}
c_G^\ast\leq d_G\leq\wdtl k_G^\ast, \tag{*}
\end{gather*}
i.e.~the M\"obius family is minimal and the Lempert family is maximal.

\halfskip

The Green function $g_G$ may be generalized as follows.
Let $\bs{p}:G\too\RR_+$ be an arbitrary function. Define
\begin{multline*}
g_G(\bs{p},z):=\sup\{u(z):\;\; u:G\too[0,1),\; \log u\in\PSH(G),\\
\forall_{a\in G}\;\exists_{C=C(u,a)>0}\;\forall_{w\in G}:\;u(w)\leq
C\|w-a\|^{\bs{p}(a)}\},\quad z\in G;\quad\footnotemark
\end{multline*}
\footnotetext{Here $0^0:=1$.}%
obviously the above growth condition is trivially satisfied at all
points $a\in G$ such that $\bs{p}(a)=0$. We have $g_G(\bs{0},\cdot)\equiv1$.
The function $g_G(\bs{p},\cdot)$ is called the {\it generalized
pluricomplex Green function with poles (weights) $\bs{p}$}.
Observe that if the set $|\bs{p}|:=\{a\in G: \bs{p}(a)>0\}$ is
not pluripolar, then $g_G(\bs{p},\cdot)\equiv0$.

In the case where the set $|\bs{p}|$ is finite, the function
$g_G(\bs{p},\cdot)$ was introduced by P.~Lelong in \cite{Lel1989}.

In the case where $\bs{p}=\chi_A=$ the characteristic function of a
set $A\subset G$, we put $g_G(A,\cdot):=g_G(\chi_A,\cdot)$. Obviously,
$g_G(\{a\},\cdot)=g_G(a,\cdot)$, $a\in G$.

The generalized Green function was recently studied by many authors,
e.g. \cite{CarWie2002}, \cite{Com2000a}, \cite{Edi2002}, \cite{EdiZwo1998b},
\cite{LarSig1998b}.

\halfskip

Using the same idea, one can generalize the M\"obius function.
For
$$
\bs{p}:G\too\ZZ_+:=\{0,1,2,\dots\}
$$
we put
$$
m_G(\bs{p},z):=\sup\{|f(z)|: f\in\CO(G,E),\; \ord_a\!f\geq\bs{p}(a),\;a\in G\},
\quad z\in G.
$$
The function $m_G(\bs{p},\cdot)$ is called the {\it generalized M\"obius
function with weights $\bs{p}$}. Clearly $m_G(\bs{0},\cdot)\equiv1$.
Observe that if the set $|\bs{p}|$ is not thin, then
$m_G(\bs{p},\cdot)\equiv0$. Similarly as in the case of the generalized
Green function we put $m_G(A,\cdot):=m_G(\chi_A,\cdot)$, $A\subset G$.
We get $m_G(\{a\},\cdot)=c_G^\ast(a,\cdot)$, $a\in G$.
Moreover, if $|\bs{p}|=\{a\}$ and $\bs{p}(a)=k$, then
$m_G(\bs{p},\cdot)=[m_G^{(k)}(a,z)]^k$.

It is clear that $m_G(\bs{p},\cdot)\leq g_G(\bs{p},\cdot)$ (for any function
$\bs{p}:G\too\ZZ_+$).
Some other properties of $g_G(\bs{p},\cdot)$ and $m_G(\bs{p},\cdot)$ will be
presented in \S\;2.

\halfskip

Consider the following definition.
\begin{Definition}\label{Def}
A family $\underline{d}=(d_G)_G$ of functions
$$
d_G:\RR_+^G\times G\too\RR_+
$$
is said to be {\it a generalized holomorphically contractible family} if the
following three axioms are satisfied:
\begin{itemize}
\item[(E)] $\prod_{a\in E}[m_E(a,z)]^{\bs{p}(a)}\leq d_E(\bs{p},z)\leq
\inf_{a\in E}[m_E(a,z)]^{\bs{p}(a)}$ for every $(\bs{p},z)\in\RR_+^ E\times E$
\footnote{We put $\prod_{a\in A}h(a):=\inf\{\prod_{a\in B}h(a): B\subset A,\;
\#B<+\infty\}$, $h:A\too[0,1]$.}\!,

\item[(H)] for any $F\in\CO(G,D)$ and $\bs{q}:D\too\RR_+$ we have
$d_D(\bs{q},F(z))\leq d_G(\bs{q}\circ F,z)$ for every $z\in G$,

\item[(M)] for any $\bs{p}, \bs{q}:G\too\RR_+$, if $\bs{p}\leq\bs{q}$, then
$d_G(\bs{q},\cdot)\leq d_G(\bs{p},\cdot)$.
\end{itemize}

If in the above definition one considers only integer valued weights
(like in the case of the generalized M\"obius function),
then we get the definition of a
{\it generalized holomorphically contractible family with integer valued
weights}.

Put $d_G(A,\cdot):=d_G(\chi_A,\cdot)$, $A\subset G$, $d_G(a,\cdot):=
d_G(\{a\},\cdot)$, $a\in G$.
\end{Definition}

One can prove that the generalized Green and M\"obius
functions satisfy all the above axioms (cf.~\S\;2).

The main result of the paper is the following theorem.

\begin{Theorem}\label{T1}
In the category of generalized holomorphically contractible families there
exists a minimal and a maximal object. They are given by the following formulae:
\begin{align*}
d_G^{\min}(\bs{p},z):=&\sup\{
\prod_{\mu\in f(G)}[m_E(\mu,f(z))]^{\sup\bs{p}(f^{-1}(\mu))}: f\in\CO(G,E)\}\\
=&\sup\{
\prod_{\mu\in f(G)}|\mu|^{\sup\bs{p}(f^{-1}(\mu))}: f\in\CO(G,E),\;f(z)=0\},\\
d_G^{\max}(\bs{p},z):=&\inf\{[\wdtl k_G^\ast(a,z)]^{\bs{p}(a)}: a\in G\}\\
=&\inf\{|\mu|^{\bs{p}(\phi(\mu))}: \phi\in\CO(E,G),\;\phi(0)=z,\;\mu\in E\}.
\end{align*}
\end{Theorem}

Observe that if $|\bs{p}|=\{a\}$ and $\bs{p}(a)=k$, then
$d_G^{\min}(\bs{p},\cdot)=[c_G^\ast(a,\cdot)]^k$ and
$d_G^{\max}(\bs{p},\cdot)=[\wdtl k_G^\ast(a,\cdot)]^k$. Moreover, for
$A\subset G$ we get
\begin{align*}
d_G^{\min}(A,z)=&\sup\{\prod_{\mu\in f(A)}m_E(\mu,f(z)): f\in\CO(G,E)\}\\
=&\sup\{\prod_{\mu\in f(A)}|\mu|: f\in\CO(G,E),\;f(z)=0\},\\
d_G^{\max}(A,z)=&\inf\{\wdtl k_G^\ast(a,z): a\in A\}.
\end{align*}
The function $d_G^{\min}$ (resp.~$d_G^{\max}$) may be considered as a
generalization of the M\"obius function $c_G^\ast$ (resp.~Lempert function
$\wdtl k_G^\ast$).
The proof of Theorem \ref{T1} will be given in \S\;3. Some properties of
$d_G^{\min}$ and $d_G^{\max}$ will be presented in \S\;4.

\halfskip

\section{Basic properties of $g_G$ and $m_G$.}

Directly from the definitions we conclude that the systems
$(g_G)_G$ and $(m_G)_G$ satisfy (H) and (M) and the following conditions
(to simplify formulations we will write
$d_G$ if a given property holds  simultaneously for $m_G$ and $g_G$):

\begin{ETh}\label{2.1}
$d_G(\bs{p},\cdot) d_G(\bs{q},\cdot)\leq d_G(\bs{p}+\bs{q},\cdot)\leq
\min\{d_G(\bs{p},\cdot),\; d_G(\bs{q},\cdot)\}$. In particular,
$g_G(\bs{p},\cdot)\leq\inf_{a\in G}[g_G(a,\cdot)]^{\bs{p}(a)}
\leq d_G^{\max}(\bs{p},\cdot)$.
\end{ETh}

\begin{ETh}\label{2.2}
If the set $|\bs{p}|$ is finite, then
$\prod_{a\in G}[d_G(a,\cdot)]^{\bs{p}(a)}\leq d_G(\bs{p},\cdot)$.
\end{ETh}

\begin{ETh}\label{2.3}
$g_G(\bs{p},z)=\sup\{u(z):\;\; u:G\too[0,1),\;\log u\in\PSH(G)$,

\hfill $u(\cdot)\leq\inf_{a\in G}[g_G(a,\cdot)]^{\bs{p}(a)}\}$, $z\in G$.
\end{ETh}

\begin{ETh}\label{2.4}
$m_G(\bs{p},\cdot)\in\Cal C(G)$.
\end{ETh}

\begin{proof}
The family $\{f\in\CO(G,E): \ord_a\!f\geq\bs{p}(a),\;a\in G\}$
is equicontinuous.
\end{proof}

\begin{ETh}\label{2.5}
If $\bs{p}\not\equiv0$, then for any $z_0\in G$ there exists an extremal function for $m_G(\bs{p},z_0)$,
i.e.~a function  $f_{z_0}\in\CO(G,E)$, $\ord_a\!f\geq\bs{p}(a)$, $a\in G$,
and $m_G(\bs{p},z_0)=|f_{z_0}(z_0)|$.
\end{ETh}

\begin{ETh}\label{2.6}
$\log d_G(\bs{p},\cdot)\in\PSH(G)$.
\end{ETh}

\begin{proof} We can argue as in the one-pole case ---
cf.~\cite{JarPfl1993}, \S\S \;2.5, 4.2.
\end{proof}

\begin{ETh}\label{2.7}
If $G_k\nearrow G$ and $\bs{p}_k\nearrow\bs{p}$, then
$d_{G_k}(\bs{p}_k,z)\searrow d_G(\bs{p},z)$, $z\in G$.
\end{ETh}

\begin{proof} It is clear that the sequence is monotone and the limit
function $u$ satisfies $u\geq d_G(\bs{p},\cdot)$.

In the case of the generalized Green function, using \ref{2.6}, we have
$u\in\PSH(G)$. By \ref{2.3} it remains to observe that
$u(z)\leq\inf_{a\in G}[g_G(a,z)]^{\bs{p}(a)}$, $z\in G$
(because $g_{G_k}(a,z)\searrow g_G(a,z)$ for every $(a,z)\in G\times G$).

The case of the generalized M\"obius function is simpler and it follows from
\ref{2.5} and a Montel argument.
\end{proof}

\begin{ETh}\label{2.8}
$g_G(\bs{p},\cdot)=\inf\{g_G(\bs{q},\cdot): \bs{q}\leq\bs{p}, \;\#|\bs{q}|<+\infty\}$.
\end{ETh}

\begin{proof}
Let $u:=\inf\{g_G(\bs{q},\cdot): \bs{q}\leq\bs{p}, \;\#|\bs{q}|<+\infty\}$.
Obviously $g_G(\bs{p},\cdot)\leq u$. By \ref{2.3}, to prove the opposite inequality we only
need to show that $\log u$ is plurisubharmonic. Observe that
$g_G(\max\{\bs{q}_1,\dots,\bs{q}_N\},\cdot)
\leq\min\{g_G(\bs{q}_1,\cdot),\dots,g_G(\bs{q}_N,\cdot)\}$.
We finish the proof by applying the following general result.

\halfskip
\begin{Lemma}\label{infsub}
Let $(v_i)_{i\in A}\subset\PSH(\Omega)$ ($\Omega\subset\CC^n$) be such that
for any $i_1,\dots,i_N\in A$ there exists an $i_0\in A$ such that
$v_{i_0}\leq\min\{v_{i_1},\dots,v_{i_N}\}$. Then $v:=\inf_{i\in A}v_i\in\PSH(\Omega)$.
\end{Lemma}

\begin{proof} It suffices to consider the case $n=1$.
Take a disc $\Delta_a(r)\Subset\Omega$, $\eps>0$, and a continuous function
$w\in\Cal C(\partial\Delta_a(r))$ such that $w\geq v$ on $\partial\Delta_a(r)$.
We want to show that $v(a)\leq\frac1{2\pi}\int_0^{2\pi}w(a+re^{i\theta})d\theta+\eps$.
For any point $b\in\partial\Delta_a(r)$ there exists an $i=i(b)\in A$ such that
$v_i(b)< w(b)+\eps$. Hence there exists an open arc $I=I(b)\subset
\partial\Delta_a(r)$ with $b\in I$ such that $v_i(\lambda)<w(\lambda)+\eps$, $\lambda\in I$.
By a compactness argument, we find $b_1,\dots,b_N\in\partial\Delta_a(r)$
such that $\partial\Delta_a(r)=\bigcup_{j=1}^NI(b_j)$. By assumption, there
exists an $i_0\in A$ such that $v_{i_0}\leq\min\{v_{i(b_1)},\dots,v_{i(b_N)}\}$.
Then
\begin{gather*}
v(a)\leq v_{i_0}(a)\leq
\frac1{2\pi}\int_0^{2\pi}v_{i_0}(a+re^{i\theta})d\theta\leq
\frac1{2\pi}\int_0^{2\pi}w(a+re^{i\theta})d\theta+\eps.
\qquad{\text{\qed}}\Dqed
\end{gather*}
\end{proof}
\end{proof}

\begin{ETh}\label{2.9}
$\prod_{a\in G}[g_G(a,\cdot)]^{\bs{p}(a)}\leq g_G(\bs{p},\cdot)$.
\end{ETh}

\begin{proof} Use \ref{2.2} and \ref{2.8}.
\end{proof}

\begin{ETh}\label{2.10}
If $G\subset\CC$, then
$
g_G(\bs{p},z)=\prod_{a\in G}[g_G(a,z)]^{\bs{p}(a)}$, $z\in G$.

In particular,
$g_E(\bs{p},z)=\prod_{a\in E}[m_E(a,z)]^{\bs{p}(a)}$, $z\in E$.
\end{ETh}

\begin{proof}
By \ref{2.8} we may assume that the set $|\bs{p}|$ is finite.
Now, by \ref{2.7}, we may assume that $G\Subset\CC$ is regular with respect to
the Dirichlet problem.
Let $u:=\prod_{a\in|\bs{p}|}[g_G(a,\cdot)]^{\bs{p}(a)}$. Then
the function $\log u$ is subharmonic on $G$ and harmonic on
$G\setminus|\bs{p}|$. The function $v:=\log g_G(\bs{p},\cdot)-\log u$
is locally bounded from above in $G$ and $\limsup_{z\to\zeta}v(z)\leq0$,
$\zeta\in\partial G$. Consequently, $v$ extends to a subharmonic function on
$G$ and, by the maximum principle, $v\leq0$ on $G$,
i.e.~$g_G(\bs{p},\cdot)\leq u$ on $G$. The opposite inequality follows from
\ref{2.9}.
\end{proof}

\begin{ETh}\label{2.11} For any function $\bs{p}: G\too\ZZ_+$ we get
$$
m_G(\bs{p},\cdot)=
\inf\{m_G(\bs{q},\cdot): \bs{q}:G\too\ZZ_+,\;
\bs{q}\leq\bs{p},\;\#|\bs{q}|<+\infty\}.
$$
In particular, for any function $\bs{p}:E\too\ZZ_+$ we have
$$
m_E(\bs{p},z)=g_E(\bs{p},z)=\prod_{a\in E}[m_E(a,z)]^{\bs{p}(a)},\quad z\in E.
$$
\end{ETh}

\begin{proof}  The case where $|\bs{p}|$ is finite is trivial.
The case where the set $|\bs{p}|$ is countable follows from \ref{2.7}.
In the general case let $A_k:=\{a\in G: \bs{p}(a)=k\}$ and let $B_k$
be a countable (or finite) dense subset of $A_k$, $k\in\ZZ_+$. Put
$B:=\bigcup_{k=0}^\infty B_k$, $\bs{p}':=\bs{p}\cdot\chi_B$. Then
$\bs{p}'\leq\bs{p}$, the set $|\bs{p}'|$ is countable, and $m_G(\bs{p},\cdot)
\equiv m_G(\bs{p}',\cdot)$. Consequently, the result reduces to the
countable case.
\end{proof}

\begin{Proposition}[\cite{EdiZwo1998b}, \cite{LarSig1998b}]\label{GGF3}
Let $G, D\subset\CC^n$ be domains
and let $F:G\too D$ be a proper holomorphic mapping. Let
$\bs{q}:D\too\RR_+$. Assume that $\det F'(a)\neq0$ for any
$a\in G$ such that $\bs{q}(F(a))>0$. Then
$$
g_D(\bs{q},F(z))=g_G(\bs{q}\circ F,z),\quad z\in G.
$$
In particular, if $B\subset D$ is such that $\det F'(a)\neq0$ for any
$a\in F^{-1}(B)$, then
$$
g_D(B,F(z))= g_G(F^{-1}(B),z),\quad z\in G.
$$
\end{Proposition}

\begin{Corollary}\label{GGF4}
Let $A_1,\dots,A_n\subset E$ be finite sets. Put
\begin{gather*}
F_j(\lambda):=\prod_{a\in A_j}\frac{\lambda-a}{1-\overline{a}\lambda},
\quad\lambda\in E,\;j=1,\dots,n,\\
F(z):=(F_1(z_1),\dots,F_n(z_n)),\quad z=(z_1,\dots,z_n)\in E^n.
\end{gather*}
Then
\begin{align*}
m_{E^n}(A_1\!\times\!\dots\!\times\! A_n,z)&\leq g_{E^n}(A_1\!\times\!\dots\!\times\! A_n,z)\\
&=g_{E^n}(0,F(z))=\max\{|F_j(z_j)|: j=1,\dots,n\}\\
&=\max\{m_E(A_1,z_1),\dots,m_E(A_n,z_n)\}\\
&\leq m_{E^n}(A_1\!\times\!\dots\!\times\! A_n,z),\quad z=(z_1,\dots,z_n)\!\in\! E^n.
\end{align*}
\end{Corollary}

\begin{Proposition}[\cite{CarWie2002}]\label{GGF9}
Let $\bs{p}:E^n\too\RR_+$ be such that
$|\bs{p}|=\{a_1,\dots,a_N\}\subset E\times\{0\}^{n-1}$. Put $k_j:=\bs{p}(a_j)$,
$j=1,\dots,N$, and assume that $k_1\geq\dots\geq k_N$. Then
$$
g_{E^n}(\bs{p},z)=\prod_{j=1}^N u_j^{k_j-k_{j+1}}(z),\quad z\in E^n,
$$
where $k_{N+1}:=0$ and
\begin{align*}
u_j(z):&=\max\{m_E(a_{1,1},z_1)\cdots m_E(a_{j,1},z_1),\;|z_2|,\dots,|z_n|\}\\
&=
\max\{m_E(\{a_{1,1},\dots,a_{j,1}\},z_1),|z_2|,\dots,|z_n|\}\\
&=g_{E^n}(\{a_1,\dots,a_j\},z),\quad j=1,\dots,N.
\end{align*}

If $k_1,\dots,k_N\in\NN$, then $m_{E^n}(\bs{p},\cdot)=g_{E^n}(\bs{p},\cdot)$.
\end{Proposition}

Observe that if $k_1=\dots=k_N=1$, then the above formula coincides
with that from Corollary \ref{GGF4}.

Notice that, even for the simplest case not covered by Proposition \ref{GGF9}
$n=N=2$, $a_1=(0,0)$, $a_2\in(E_\ast)^2$, $k_1=k_2=1$, an effective
formula for $g_{E^n}(\bs{p},\cdot)$ is not known.

\halfskip

Recall that by the Lempert theorem (cf.~\cite{JarPfl1993}, Ch.~8), if
$G\subset\CC^n$ is convex, then $c_G^\ast=\wdtl k_G^\ast$, and consequently,
by (*), all holomorphically contractible families coincide on $G$. The following
example shows that this is not true in the category of generalized holomorphically
contractible families.

\begin{Example}[Due to W.~Zwonek]\label{Zwonek}
Let $D:=\{(z,w)\in\CC^2: |z|+|w|<1\}$, $A_t:=\{(t,\sqrt{t}), (t,-\sqrt{t})\}$,
$0<t\ll1$. Then
$$
m_D(A_t,(0,0))<g_D(A_t,(0,0))<d_D^{\max}(A_t,(0,0))
$$
for small $t$.

\halfskip

Indeed, let $G:=\{(z,w)\in\CC^2: |z|+\sqrt{|w|}<1\}$ and let $F:D\too G$,
$F(z,w):=(z,w^2)$. Note that $F$ is proper and locally biholomorphic
in a neighborhood of $A_t$. Moreover, $A_t=F^{-1}(t,t)$.

Using Proposition \ref{GGF3}, we conclude
that $g_D(A_t,(0,0))=g_G((t,t),(0,0))$.

Observe that $m_D(A_t,(0,0))=m_G((t,t),(0,0))$. In fact, the inequality
`$\geq$' follows from (H) (applied to $F$). The opposite inequality may be
proved as follows. Let $f\in\CO(D,E)$ be such that $f|_{A_t}=0$. Define
$$
\wdtl f(z,w):=(1/2)(f(z,\sqrt{w})+f(z,-\sqrt{w})),\quad (z,w)\in G.
$$
Note that $\wdtl f$ is well defined, $|\wdtl f|<1$, $\wdtl f(t,t)=0$,
$\wdtl f$ is continuous,
and $\wdtl f$ is holomorphic on $D\cap\{w\neq0\}$.
In particular, $\wdtl f$ is holomorphic on $D$.
Consequently, $|f(0,0)|=|\wdtl f(0,0)|\leq m_G((t,t),(0,0))$.

Suppose that $m_D(A_{t_k},(0,0))=g_D(A_{t_k},(0,0))$ for a sequence $t_k\searrow0$.
Then
\begin{multline*}
g_G((t_k,t_k),(0,0))=g_D(A_{t_k},(0,0))=m_D(A_{t_k},(0,0))\\
=m_G((t_k,t_k),(0,0))\leq g_G((t_k,t_k),(0,0)),\quad k=1,2,\dots.
\end{multline*}
Thus $m_G((t_k,t_k),(0,0))=g_G((t_k,t_k),(0,0))$, $k=1,2,\dots$.

Consequently, using \cite{JarPfl1993}, \S\;2.5, and \cite{Zwo2000c}
(Corollary 4.4) (or \cite{Zwo2000c}, Corollary 4.2.3), we conclude that
$$
\gamma_G((0,0);(1,1))=A_G((0,0);(1,1)),
$$
where $\gamma_G$ (resp.~$A_G$) denotes the Carath\'eodory--Reiffen
(resp.~Azukawa) metric of
$G$ (cf.~\cite{JarPfl1993}, \S\S\;2.1, 4.2). Hence, by
Propositions 4.2.7 and 2.2.1(d) from \cite{JarPfl1993}, using the fact that
$D$ is the convex envelope of $G$, we get
\begin{gather*}
2=h_D(1,1)
=\gamma_G((0,0);(1,1))=A_G((0,0);(1,1))=h_G(1,1)=\frac{2}{3-\sqrt{5}},
\end{gather*}
where $h_D$ (resp.~$h_G$) denotes the Minkowski function for $D$ (resp.~$G$);
contradiction.

\halfskip

To see the inequality $g_D(A_t,(0,0))<d_D^{\max}(A_t,(0,0))$, we may argue
as follows.

We already know that
$$
g_D(A_t,(0,0))=g_G((t,t),(0,0))\approx
g_G((0,0),(t,t))=h_G(t,t)=\frac{2t}{3-\sqrt{5}},\quad t\approx0.
$$

On the other hand
\begin{multline*}
d_D^{\max}(A_t,(0,0))=
\min\{\wdtl k_D^\ast((t,-\sqrt{t}), (0,0)), \wdtl k_D^\ast((t,\sqrt{t}),(0,0))\}\\
=\min\{h_D(t,-\sqrt{t}), h_D(t,\sqrt{t})\}=t+\sqrt{t}.
\end{multline*}

It remains to observe that $\frac{2t}{3-\sqrt{5}}< t+\sqrt{t}$ for small $t>0$.

\halfskip

Let $\delta_D(A_t,\cdot)$ denote the
{\it Coman function for $D$ with poles at $A_t$}, i.e.
\begin{multline*}
\delta_D(A_t,(z,w))=\inf\{|\mu_1\mu_2|: \exists_{\phi\in\CO(E,D)}:\;\\
\phi(0)=(z,w),\;\phi(\mu_1)=(t,\sqrt{t}),\;\phi(\mu_2)=(t,-\sqrt{t})\},\quad (z,w)\in D,
\end{multline*}
cf.~\cite{Com2000a}. It is known that $g_D(A_t,\cdot)\leq\delta_D(A_t,\cdot)$.
Taking $\phi(\lambda):=(\lambda^2/4,\lambda/2)$, we easily see that
$\delta_D(A_t,(0,0))\leq4t<t+\sqrt{t}=d_D^{\max}(A_t,(0,0))$, $0<t\ll1$.
{\it We do not know whether $g_D(A_t,(0,0))<\delta_D(A_t,(0,0))$ for small
$t>0$.}
\end{Example}

\halfskip

\section{Proof of Theorem 1.2.}

{\bf Step 1.} {\it If $(d_G)_G$ satisfies {\rm (H)} and
\begin{itemize}
\item[{\rm (E${}^+$)}] $d_E(\bs{p},\lambda)\leq d_E^{\max}(\bs{p},\lambda)=
\inf\{[m_E(\mu,\lambda)]^{\bs{p}(\mu)}: \mu\in E\},\quad(\bs{p},\lambda)\in
\RR_+^E\times E$,
\end{itemize}
\noindent then $d_G\leq d_G^{\max}$ for any $G$.

The result remains true in the category of contractible families with integer
valued weights.}
\begin{proof}
\begin{align*}
d_G(\bs{p},z)\overset{\text{\rm (H)}}\leq&\inf\{d_E(\bs{p}\circ\phi,0):
\phi\in\CO(E,G),\;\phi(0)=z\}\\
\overset{\text{\rm (E${}^+$)}}
\leq&\inf\{|\mu|^{\bs{p}(\phi(\mu))}: \phi\in\CO(E,G),\;\phi(0)=z,\;
\mu\in E\}\\
=&\; d_G^{\max}(\bs{p},z),\quad(\bs{p},z)\in\RR_+^G\times G.\Dqed
\end{align*}
\end{proof}

{\bf Step 2.} {\it The system $(d_G^{\max})_G$ satisfies {\rm (E)},
{\rm (H)}, and {\rm (M)}.}
\begin{proof} (E) and (M) are obvious. To prove (H) let $F:G\too D$ be holomorphic and let
$\bs{q}:D\too\RR_+$. Then
\begin{align*}
d_D^{\max}(\bs{q},F(z))&=\inf\{[\wdtl k_D^\ast(b,F(z))]^{\bs{q}(b)}: b\in D\}\\
&\leq\inf\{[\wdtl k_D^\ast(F(a),F(z))]^{\bs{q}(F(a))}: a\in G\}\\
&\leq\inf\{[\wdtl k_G^\ast(a,z)]^{\bs{q}(F(a))}: a\in G\}
=d_G^{\max}(\bs{q}\circ F,z),\quad z\in G.\Dqed
\end{align*}
\end{proof}

{\bf Step 3.} {\it If $(d_G)_G$ satisfies {\rm (H)}, {\rm (M)}, and
\begin{itemize}
\item[{\rm (E${}^-$)}] $\prod_{\mu\in E}
[m_E(\mu,\lambda)]^{\bs{p}(\mu)}\leq d_E(\bs{p},\lambda),
\quad(\bs{p},\lambda)\in \RR_+^E\times E$,
\end{itemize}
\noindent then $d_G^{\min}\leq d_G$ for any $G$.

The result remains true in the category of contractible families with integer
valued weights.}
\begin{proof}
\begin{align*}
d_G(\bs{p},z)\overset{\text{\rm (M)}}\geq
&\sup\{d_G(\bs{q}\circ f,z): f\in\CO(G,E),\;\bs{q}:E\too\RR_+,\;f(z)=0,\;\bs{p}\leq\bs{q}\circ f\}\\
\overset{\text{\rm (H)}}\geq
&\sup\{d_E(\bs{q},0): f\in\CO(G,E),\;\bs{q}:E\too\RR_+,\;f(z)=0,\;\bs{p}\leq\bs{q}\circ f\}\\
\overset{\text{\rm (E${}^-$)}}
\geq&\sup\{\prod_{\mu\in E}|\mu|^{\bs{q}(\mu)}:f\in\CO(G,E),\;\bs{q}:E\too\RR_+,\;f(z)=0,\;\bs{p}\leq\bs{q}\circ f\}\\
\geq&\sup\{\prod_{\mu\in f(G)}|\mu|^{\sup\bs{p}(f^{-1}(\mu))}:
f\in\CO(G,E),\;f(z)=0\}\\
=\; & d_G^{\min}(\bs{p},z),\quad (\bs{p},z)\in\RR_+^G\times G.\Dqed
\end{align*}
\end{proof}

{\bf Step 4.} {\it The system $(d_G^{\min})_G$ satisfies {\rm (E)},
{\rm (H)}, and {\rm (M)}.}

\begin{proof} (E) and (M) are elementary.
To prove (H) let $F:G\too D$ be holomorphic and let $\bs{q}:D\too\RR_+$. Then
\begin{align*}
d_D^{\min}(\bs{q},F(z))=&\sup\{\prod_{\mu\in g(D)}[m_E(\mu,g(F(z))]^{\sup\bs{q}(g^{-1}(\mu))}:
g\in\CO(D,E)\}\\
\overset{f=g\circ F}\leq
&\sup\{\prod_{\mu\in f(G)}[m_E(\mu,f(z))]^{\sup(\bs{q}\circ F)(f^{-1}(\mu))}:
f\in\CO(G,E)\}\\
=\;&d_G^{\min}(\bs{q}\circ F,z),\quad z\in G.\Dqed
\end{align*}
\end{proof}

\begin{Corollary}\label{3.1}
{\rm (a)} $d_G^{\min}\leq g_G\leq d_G^{\max}$ and
$d_G^{\min}\leq m_G\leq g_G\leq d_G^{\max}$ (for integer valued weights).

{\rm (b)} $d_E^{\min}(\bs{p},\lambda)=g_E(\bs{p},\lambda)=\prod_{\mu\in E}
[m_E(\mu,\lambda)]^{\bs{p}(\mu)},\quad(\bs{p},\lambda)\in\RR_+^E\times E$.

{\rm (c)} $d_G^{\min}(A,\cdot)=m_G(A,\cdot)$ for any $A\subset G$.
\end{Corollary}

\begin{proof} (a) follows from Theorem \ref{T1}.

(b) Using (a) and \ref{2.10} we get
$$
\prod_{\mu\in E}[m_E(\mu,\lambda)]^{\bs{p}(\mu)}\leq d_E^{\min}(\bs{p},\lambda)
\leq g_G(\bs{p},\lambda)=\prod_{\mu\in E}[m_E(\mu,\lambda)]^{\bs{p}(\mu)}.
$$

(c) Let $A\subset G$. Then
\begin{align*}
m_G(A,z)&\geq d_G^{\min}(A,z)
\geq\sup\{\prod_{\mu\in f(A)}m_E(\mu,f(z)): f\in\CO(G,E),\;f|_A=0\}\\
&=m_G(A,z),\quad z\in G.\Dqed
\end{align*}
\end{proof}

\begin{Example}
Let $G:=E^2$, $|\bs{p}|=\{(-\frac12,0), (\tfrac12,0)\}$, $\bs{p}(-\frac12,0)=2$,
$\bs{p}(\frac12,0)=1$. Then
$d_{E^2}^{\min}(\bs{p},(0,\frac13))<m_{E^2}(\bs{p},(0,\frac13))$
(cf.~Corollary \ref{3.1}(c)).

\halfskip

Indeed, by Proposition \ref{GGF9},
$$
m_{E^2}((\bs{p},(0,\tfrac13))=u_1(0,\tfrac13)u_2(0,\tfrac13)=
\max\{\tfrac12,\tfrac13\}\max\{\tfrac12\cdot\tfrac12,\tfrac13\}=
\tfrac12\cdot\tfrac13=\tfrac16.
$$
On the other side:
\begin{align*}
&d_{E^2}^{\min}(\bs{p},(0,\tfrac13))\\
&=\max\{\sup\{|f(-\tfrac12,0)|^2|f(\tfrac12,0)|:
f\in\CO(E^2,E), f(0,\tfrac13)=0, f(-\tfrac12,0)\neq f(\tfrac12,0)\},\\
&\;\;\quad\qquad\sup\{|f(0,\tfrac13)|^2: f\in\CO(E^2,E),\;f(-\tfrac12,0)=f(\tfrac12,0)=0\}\}\\
&\leq\max\{[m_{E^2}((-\tfrac12,0),(0,\tfrac13))]^2
m_{E^2}((\tfrac12,0),(0,\tfrac13))\},\\
&\;\;\quad\qquad[m_{E^2}(\{(-\tfrac12,0),(\tfrac12,0)\},(0,\tfrac13))]^2\}\\
&=\max\{[\max\{\tfrac12,\tfrac13\}]^2\max\{\tfrac12,\tfrac13\},
[m_{E^2}(\{-\tfrac12,\tfrac12\}\times\{0\},(0,\tfrac13))]^2\}\\
&=\max\{\tfrac18, [\max\{\tfrac12\cdot\tfrac12,\tfrac13\}]^2\}=\tfrac18.
\end{align*}
\end{Example}

\halfskip

\section{Basic properties of $d_G^{\min}$ and $d_G^{\max}$.}

\begin{ETh}\label{4.3} If $D\subset\CC^m$ is a Liouville domain, then
$$
d_{G\times D}^{\min}(\bs{p},(z,w))=d_G^{\min}(\bs{p}',z),\quad (z,w)\in G\times D,
$$
where $\bs{p}'(z):=\sup\{\bs{p}(z,w): w\in D\}$, $z\in G$, and
$d_G^{\min}(\bs{p}',\cdot):=0$ if there exists a $z_0\in G$ with
$\bs{p}'(z_0)=+\infty$.
\end{ETh}

\begin{ETh}\label{4.4}
{\rm (a)} The functions $d_G^{\min}(\bs{p},\cdot)$ and $d_G^{\max}(\bs{p},\cdot)$
are upper semicontinuous.

{\rm (b)} If $\bs{p}:G\too\ZZ_+$, then
$d_G^{\min}(\bs{p},\cdot)\in\Cal C(G)$.
\end{ETh}

\begin{proof} (a) The case of $d_G^{\max}(\bs{p},\cdot)$ is obvious. To prove
the upper semicontinuity of $d_G^{\min}(\bs{p},\cdot)$, fix a $z_0\in G$ and
suppose that $d_G^{\min}(\bs{p},z_k)
\too\alpha>\beta>d_G^{\min}(\bs{p},z_0)$ for some sequence $z_k\too z_0$.
Let $f_k\in\CO(G,E)$ be such that $f_k(z_k)=0$ and
$\prod_{\mu\in f_k(G)}|\mu|^{\sup\bs{p}(f_k^{-1}(\mu))}\too\alpha$.
By a Montel argument we may assume that $f_k\too f_0$ locally uniformly
in $G$ with $f_0\in\CO(G,E)$, $f_0(z_0)=0$. Since
$\prod_{\mu\in f_0(G)}|\mu|^{\sup\bs{p}(f_0^{-1}(\mu))}<\beta$, we can find
a finite set $A\subset G$ such that $f_0|_A$ is injective and
$\prod_{a\in A}|f_0(a)|^{\bs{p}(a)}<\beta$. Consequently,
$\prod_{a\in A}|f_k(a)|^{\bs{p}(a)}<\beta$ and $f_k|_A$ is injective
for $k\gg1$. Finally,
$\prod_{\mu\in f_k(G)}|\mu|^{\sup\bs{p}(f_k^{-1}(\mu))}<\beta$, $k\gg1$;
contradiction.

(b) In view of (a), it suffices to prove that for every $f\in\CO(G,E)$ the
function $u_f(z):=\prod_{\mu\in f(G)}[m_E(\mu,f(z))]^{\sup\bs{p}(f^{-1}(\mu))}$, $z\in G$,
is continuous on $G$. Observe that
\begin{gather*}
u_f(z)=\inf_M\{\prod_{\mu\in M}[m_E(\mu,f(z))]^{k_f(\mu)}\},
\end{gather*}
where $M$ runs over all finite sets $M\subset f(|\bs{p}|)$ such that
$k_f(\mu):=\sup\bs{p}(f^{-1}(\mu))<+\infty$, $\mu\in M$.
Thus $u_f=\inf_M\{|h_M|\}$, where $h_M\in\CO(G,E)$. Consequently,
since the family $(h_M)_M$ is equicontinuous, the function $u_f$ is
continuous on $G$.
\end{proof}

\begin{Example}
Let $\bs{p}:E\times\CC\too\RR_+$, $\bs{p}(\frac1k,k):=\frac1{k^2}$, $k=2,3,\dots$,
and $\bs{p}(z,w):=0$ otherwise. Notice that $|\bs{p}|$ is discrete.
Then by \ref{4.3} and  Corollary \ref{3.1}(b),
$$
d_{E\times\CC}^{\min}(\bs{p},(z,w))=d_E^{\min}(\bs{p}',z)=\prod_{k=2}^\infty
[m_E(1/k,z)]^{1/k^2},\quad (z,w)\in E\times\CC.
$$
In particular, $d_{E\times\CC}^{\min}(\bs{p},\cdot)$ is discontinuous
at $(0,w)\in E\times\CC\setminus|\bs{p}|$.
\end{Example}

\begin{ETh}[Cf.~\ref{2.5}]\label{4.5}
If $\#|\bs{p}|<+\infty$, then for any $z_0\in G$ there exists an extremal function
for $d_G^{\min}(\bs{p},z_0)$, i.e.~a function $f_{z_0}\in\CO(G,E)$ with
$f_{z_0}(z_0)=0$ and
$$
\prod_{\mu\in f_{z_0}(G)}|\mu|^{\sup\bs{p}(f_{z_0}^{-1}(\mu)
)}=d_G^{\min}(\bs{p},z_0).
$$
\end{ETh}

\begin{proof} Fix a $z_0\in G$ and let $f_k\in\CO(G,E)$, $f_k(z_0)=0$ be such that
$$
\alpha_k:=\prod_{\mu\in f_k(G)}|\mu|^{\sup\bs{p}(f_k^{-1}(\mu))}\too\alpha:=d_G^{\min}(\bs{p},z_0).
$$
Let $A_k\subset|\bs{p}|$ be such that $f_k|_{A_k}$ is injective,
$f_k(A_k)=f_k(|\bs{p}|)$, and $\bs{p}(a)=\sup\bs{p}(f_k^{-1}(f_k(a)))$, $a\in A_k$.
Thus $\alpha_k=\prod_{a\in A_k}|f_k(a)|^{\bs{p}(a)}$.
We may assume that $A_k=B$ is independent of $k$ and for any
$a\in B$ the fiber $B_a:=f_k^{-1}(f_k(a))\cap|\bs{p}|$ is also independent of $k$.
Moreover, we may assume that $f_k\too f_0$
locally uniformly in $G$. Then $f_0\in\CO(G,E)$, $f_0(z_0)=0$, and
$\prod_{a\in B}|f_0(a)|^{\bs{p}(a)}=\alpha$.
Observe that $f_0(B)=f_0(|\bs{p}|)$. Let $B_0\subset B$ be such that
$f_0|_{B_0}$ is injective and $f_0(B_0)=f_0(B)$. We have
\begin{align*}
\alpha&\geq\prod_{\mu\in f_0(|\bs{p}|)}|\mu|^{\sup\bs{p}(f_0^{-1}(\mu))}
=\prod_{\mu\in f_0(B_0)}|\mu|^{\sup\bs{p}(f_0^{-1}(\mu))}\\
&=\prod_{a\in B_0}|f_0(a)|^{\max\{\bs{p}(b): b\in B,\;f_0(b)=f_0(a)\}}
\geq\prod_{a\in B}|f_0(a)|^{\bs{p}(a)}=\alpha.\Dqed
\end{align*}
\end{proof}

\begin{ETh}\label{4.6}
$\log d_G^{\min}(\bs{p},\cdot)\in\PSH(G)$.
\end{ETh}

\begin{proof} By virtue of \ref{4.4}(a), we only need to show that
for any $f\in\CO(G,E)$ the function
$u_f(z):=\prod_{\mu\in f(G)}[m_E(\mu,f(z))]^{\sup\bs{p}(f^{-1}(\mu))}$, $z\in G$,
is log--plurisubharmonic on $G$. The proof of \ref{4.4} shows that
$u_f=\inf_Mv_M$, where $v_M$ is a log-pluri\-sub\-har\-mo\-nic function given by
the formula
$v_M:=\prod_{\mu\in M}[m_E(\mu,f(z))]^{k_f(\mu)}$ and $M$ runs over a family
of finite sets as in the proof of \ref{4.4}. Observe that
$v_{M_1\cup M_2}\leq\min\{v_{M_1},\; v_{M_2}\}$. It remains to apply
Lemma \ref{infsub}.
\end{proof}

\begin{ETh}\label{4.7}
If $G_k\nearrow G$ and $\bs{p}_k\nearrow\bs{p}$, then
$$
d_{G_k}^{\min}(\bs{p}_k,z)\searrow d_G^{\min}(\bs{p},z),\quad
d_{G_k}^{\max}(\bs{p}_k,z)\searrow d_G^{\max}(\bs{p},z),\quad z\in G.
$$
\end{ETh}

\begin{proof} By (H) and (M) the sequence is monotone and for the limit
function $u$ we have $u\geq d_G^{\min}(\bs{p},\cdot)$
(resp.~$u\geq d_G^{\max}(\bs{p},\cdot)$). Fix a $z_0\in G$.

In the case of the minimal family suppose that
$u(z_0)>\alpha>d_G^{\min}(G,z_0)$.
Let $f_k\in\CO(G_k,E)$ be such that $f_k(z_0)=0$ and
$\prod_{\mu\in f_k(G_k)}|\mu|^{\sup\bs{p}_k(f_k^{-1}(\mu))}\too u(z_0)$.
By a Montel argument we may assume that $f_k\too f_0$ locally uniformly
in $G$ with $f_0\in\CO(G,E)$, $f_0(z_0)=0$. Since
$\prod_{\mu\in f_0(G)}|\mu|^{\sup\bs{p}(f_0^{-1}(\mu))}<\alpha$, we can find
a finite set $A\subset G$ such that $f|_A$ is injective and
$\prod_{a\in A}|f_0(a)|^{\bs{p}(a)}<\alpha$. Consequently,
$\prod_{a\in A}|f_k(a)|^{\bs{p}_k(a)}<\alpha$ and $f_k|_A$ in injective
for $k\gg1$. Finally,
$\prod_{\mu\in f_k(G_k)}|\mu|^{\sup\bs{p}_k(f_k^{-1}(\mu))}<\alpha$, $k\gg1$;
contradiction.

\halfskip

In the case of the maximal family for any $a\in G$ and $\eps>0$ there exists
a $k(a,\eps)\in\NN$ such that $z_0,a\in G_k$,
$\wdtl k_{G_k}^\ast(a,z_0)\leq\wdtl k_G^\ast(a,z_0)+\eps$, and
$\bs{p}_k(a)\geq\bs{p}(a)-\eps$ for $k\geq k(a,\eps)$. Hence
\begin{align*}
\inf_{k\in\NN}&d_{G_k}^{\max}(\bs{p}_k,z_0)=\inf_{k\in\NN: a\in G_k}
[\wdtl k_{G_k}^\ast(a,z_0)]^{\bs{p}_k(a)}\\
&\leq\inf_{a\in G}\inf\{[\wdtl k_G^\ast(a,z_0)+\eps]^{\bs{p}_k(a)}:
0<\eps\ll1,\; k\geq k(a,\eps)\}\\
&\leq\inf_{a\in G}\inf\{[\wdtl k_G^\ast(a,z_0)+\eps]^{\bs{p}(a)-\eps}:
0<\eps\ll1\}=d_G^{\max}(\bs{p},z_0).\Dqed
\end{align*}
\end{proof}

\begin{Example}
Let $G:=\{z\in\CC^n: |z^\alpha|<1\}$, where $\alpha=(\alpha_1,\dots,\alpha_n)
\in\NN^n$ is such that $\alpha_1,\dots,\alpha_n$ are  relatively prime. Then
$$
d_G^{\min}(\bs{p},z)=d_E^{\min}(\bs{p}',z^\alpha)=\prod_{\mu\in E}
[m_E(\mu,z^\alpha)]^{\bs{p}'(\mu)},\quad z\in G.
$$
where $\bs{p}'(\lambda)=\sup\{\bs{p}(a): a^\alpha=\lambda\}$, $\lambda\in E$,
and $d_E^{\min}(\bs{p}',\cdot):=0$ if there exists a $\lambda_0\in E$ with
$\bs{p}'(\lambda_0)=+\infty$.

Indeed, it is known that any function $f\in\CO(G,E)$ has the form
$f=g\circ\Phi$, where $\Phi(z):=z^\alpha$ and $g\in\CO(E,E)$ --- cf.~\cite{JarPfl1993},
\S\;4.4. Thus
\begin{multline*}
d_G^{\min}(\bs{p},z)=\sup\{\prod_{\mu\in g(\Phi(G))}
[m_E(\mu,g(\Phi(z)))]^{\sup\bs{p}(\Phi^{-1}(g^{-1}(\mu)))}: g\in\CO(E,E)\}\\
=\sup\{\prod_{\mu\in g(E)}
[m_E(\mu,g(\Phi(z)))]^{\sup\bs{p}'(g^{-1}(\mu))}: g\in\CO(E,E)\}=
d_E^{\min}(\bs{p}',\Phi(z)).
\end{multline*}
\end{Example}

\halfskip

\section{Product property.}

Let $\underline{d}=(d_G)_G$ be a generalized holomorphically contractible family with integer
valued weights. We say that $\underline{d}$ has the {\it product property} if
\begin{gather*}
d_{G\times D}(A\times B, (z,w))=\max\{d_G(A,z),\;d_D(B,w)\},\quad (z,w)\in G\times D,
\tag{P}
\end{gather*}
for any domains $G\subset\CC^n$, $D\subset\CC^m$ and for any sets
$\varnothing\neq A\subset G$, $\varnothing\neq B\subset D$.
Notice that the inequality `$\geq$' follows from (H) applied to
the projections $G\times D\too G$, $G\times D\too D$.
The definition applies to the standard holomorphically
contractible families and means that
$$
d_{G\times D}((a,b), (z,w))=\max\{d_G(a,z),\;d_D(b,w)\},\quad (a,b),\;
(z,w)\in G\times D.
$$
It is well known that the families $(\wdtl k_G^\ast)_G$, $(c_G^\ast)_G$,
$(g_G)_G$ have the product property --- cf.~\cite{JarPfl1993}, Ch.9,
\cite{Edi1997b}, \cite{Edi1999}, \cite{Edi2000b}.

Moreover, it is known that the higher order M\"obius functions $(m_G^{(k)})_G$
with $k\geq2$ have no product property --- cf.~\cite{JarPfl1993}, Ch.9.

Thus {\it it is natural to ask whether the minimal and maximal families
have the product property}.

\begin{Proposition} The system $(d_G^{\max})_G$ has the product property.
\end{Proposition}

\begin{proof}
Fix $(z_0,w_0)\in G\times D$ and $\eps>0$. Let $(a,b)\in A\times B$ be such that
$\wdtl k_G^\ast(a,z_0)\leq d_G^{\max}(A,z_0)+\eps$,
$\wdtl k_D^\ast(b,w_0)\leq d_G^{\max}(B,w_0)+\eps$. Then using the product
property for $(\wdtl k_G^\ast)_G$, we get
\begin{align*}
d_{G\times D}^{\max}(A\times B,(z_0,w_0))&\leq
\wdtl k_{G\times D}^\ast((a,b),(z_0,w_0))\\
&=\max\{\wdtl k_G^\ast(a,z_0),\;\wdtl k_D^\ast(b,w_0)\}\\
&\leq\max\{d_G^{\max}(A,z_0),\;d_D^{\max}(B,w_0)\}+\eps.\Dqed
\end{align*}
\end{proof}

We do not know whether the system $(d_G^{\min})_G$ has the product property.
So far we were able to manage only the case where $\#B=1$ --- see Proposition
\ref{onepoint}.
Recall that $d_G^{\min}(A,\cdot)=m_G(A,\cdot)$ --- Corollary \ref{3.1}(c).

\begin{Proposition}\label{reduction}
Assume that for any $n\in\NN$,
the system $(m_G)_G$ has the following special product property:
\begin{gather*}
|\Psi(z,w)|\leq(\max_{G\times D}|\Psi|)\max\{m_G(A,z),\;m_D(B,w)\},\quad (z,w)\in G\times D,
\tag{P${}_0$}
\end{gather*}
where $G,D\subset\CC^n$ are balls with respect to arbitrary $\CC$--norms,
$A\subset D$, $B\subset G$ are finite and non-empty, $\Psi(z,w):=\sum_{j=1}^n
z_jw_j$, and $\Psi|_{A\times B}=0$.
Then the system $(m_G)_G$ has the product property {\rm (P)} in the full
generality.

Moreover, if {\rm (P${}_0$)} holds with $\#B=1$, then {\rm (P)} holds with
$\#B=1$.
\end{Proposition}

\begin{proof} (Cf. \cite{JarPfl1993}, the proof of Th.~9.5.)
Fix arbitrary domains $G\subset\CC^n$, $D\subset\CC^m$, non-empty sets $A\subset G$,
$B\subset G$, and $(z_0,w_0)\in G\times D$.
We have to prove that for any $F \in \CO(G\times D, E)$ with
$F|_{A\times B}=0$ the following inequality is true:
$$
|F(z_0,w_0)| \leq \max \{ m_G(A, z_0), m_D(B,w_0)\}.
$$
By \ref{2.11}, we may assume that $A, B$ are finite.

Let $(G_\nu)^\infty_{\nu=1}$, $(D_\nu)_{\nu=1}^\infty$ be sequences of
relatively compact subdomains of $G$ and $D$, respectively,
such that $A\cup\{z_0\}\subset G_\nu \nearrow G$,
$B\cup\{w_0\}\subset D_\nu \nearrow D$. By \ref{2.7}, it suffices
to show that
$$
|F(z_0,w_0)|\leq\max\{m_{G_\nu}(A,z_0), m_{D_\nu}(B,w_0)\},
\quad\nu\geq1.
$$

Fix a $\nu_0\in\NN$ and let $G':=G_{\nu_0}$, $D':=D_{\nu_0}$.

It is well known that $F$ may be approximated locally uniformly in
$G\times D$ by functions of the form
\begin{gather*}
F_s (z,w) = \sum^{N_s}_{\mu=1} f_{s,\mu} (z) g_{s,\mu} (w),
\quad  (z,w) \in G\times D,\tag{**}
\end{gather*}
where $f_{s,\mu}\in\CO(G)$, $g_{s,\mu}\in\CO(D)$, $s\geq 1$,
$\mu=1,\dots, N_s$. Notice that $F_s\too0$ uniformly on $A\times B$.
Using Lagrange interpolation formula, we find polynomials $P_s:\CC^n\times
\CC^m\too\CC$ such that $P_s|_{A\times B}=F_s|_{A\times B}$ and $P_s\too0$
locally uniformly in $\CC^n\times\CC^m$.
The functions $\wdht F_s:= F_s-P_s$, $s \geq 1$, also have the form (**) and
$\wdht F_s\too F$ locally uniformly in $G\times D$. Hence, without loss of
generality, we may assume that $F_s|_{A\times B}=0$, $s\geq 1$. Let
$m_s :=\max\{ 1,\| F_s\|_{G'\times D'}\}$ and $\wdtl F_s := F_s /m_s$,
$s\geq 1$. Note that $m_s\too 1$, and therefore
$\wdtl F_s\too F$ uniformly on $G'\times D'$. Consequently, we may
assume that $F_s (G'\times D')\Subset E$, $s\geq 1$.

It is enough to prove that
$$
|F_s(z_0, w_0)|\leq\max\{m_{G'}(A,z_0), m_{D'}(B,w_0)\},\quad s\geq1.
$$

Fix an $s=s_0\in\NN$ and let $N:=N_{s_0}$, $f_\mu:=f_{s_0,\mu}$,
$g_\mu:=g_{s_0,\mu}$, $\mu=1,\dots,N$. Let
$f:=(f_1, \dots, f_N):G\too\CC^N$ and $g:=(g_1,\dots,g_N):D\too\CC^N$.
Put
\begin{multline*}
K:=\{\xi=(\xi_1,\dots,\xi_N)\in\CC^N:\\
 |\xi_\mu|\leq\|f_\mu\|_{G'},\;
\mu=1,\dots, N,\; |\Psi(\xi,g(w))|\leq 1,\; w\in D'\}.
\end{multline*}
It is clear that $K$ is an absolutely convex compact subset of
$\CC^N$ with $f(G')\subset K$. Let
\begin{multline*}
L:=\{\eta =(\eta_1,\dots,\eta_N)\in\CC^N:\\
 |\eta_\mu|\leq\| g_\mu\|_{D'}, \; \mu = 1, \dots, N, \;
|\Psi (\xi,\eta)| \leq 1, \; \xi\in K\}.
\end{multline*}
Then again $L$ is an absolutely convex compact subset of $\CC^N$, and
moreover, $g(D')\subset L$.

Let $(W_\sigma)^\infty_{\sigma=1}$
(resp.~$(V_\sigma)^\infty_{\sigma=1}$) be a sequence of
absolutely convex bounded domains in $\CC^N$ such that
$W_{\sigma+1}\Subset W_\sigma$ and $W_\sigma\searrow K$
(resp.~$V_{\sigma + 1}\Subset V_\sigma$ and
$V_\sigma\searrow L$). Put
$M_\sigma := \|\Psi\|_{W_\sigma\times V_\sigma}$, $\sigma\in\NN$.
By (P${}_0$) and by the holomorphic contractibility applied to the mappings
$f:G'\too W_\sigma$, $g:D'\too V_\sigma$ we have
\begin{align*}
|F_{s_0}(z_0,w_0)|&=|\Psi(f(z_0),g(w_0))|\\
&\leq M_\sigma\max\{m_{W_\sigma}(f(A),f(z_0)),m_{V\sigma}(g(B),g(w_0))\}\\
&\leq M_\sigma\max\{m_{G'}(f^{-1}(f(A)),z_0),m_{D'}(g^{-1}(g(B)),w_0)\}\\
&\leq M_\sigma\max\{m_{G'}(A,z_0),m_{D'}(B,w_0)\}.
\end{align*}
Letting $\sigma\too+\infty$ we get the required result.
\end{proof}

\begin{Proposition}\label{onepoint}
The system $(m_G)_G$ has the product property {\rm (P)}
whenever $\#B=1$, i.e. for any domains $G\subset\CC^n$, $D\subset\CC^m$, for any
set $A\subset G$, and for any point $b\in D$ we have
$$
m_{G\times D}(A\times\{b\},(z,w))=\max\{m_G(A,z),\;m_D(b,w)\},
\quad (z,w)\in G\times D.
$$
\end{Proposition}

\begin{proof} By Proposition \ref{reduction}, we only
need to check (P) in the case, where $D$ is a bounded convex domain, $A$
is finite, and $B=\{b\}$. Fix $(z_0,w_0)\in G\times D$.
Let $\phi:E\too D$ be a holomorphic
mapping such that $\phi(0)=b$ and $\phi(m_D(b,w_0))=w_0$ (cf.~\cite{JarPfl1993},
Ch.~8). Consider the mapping
$F:G\times E\too G\times D$, $F(z,\lambda):=(z,\phi(\lambda))$. Then
$$
m_{G\times D}(A\times\{b\},(z_0,w_0))\leq m_{G\times E}(A\times\{0\},
(z_0,m_G(b,w_0))).
$$
Consequently, it suffices to show that
\begin{gather*}
m_{G\times E}(A\times\{0\},(z_0,\lambda))\leq\max\{m_G(A,z_0),|\lambda|\},\quad
\lambda\in E.\tag{\dag}
\end{gather*}
The case where $m_G(A,z_0)=0$ is elementary: for an $f\in\CO(G\times E,E)$
with $f|_{A\times\{0\}}=0$ we have $f(z_0,0)=0$ and hence $|f(z_0,\lambda)|\leq
|\lambda|$, $\lambda\in E$ (by the Schwarz lemma). Thus, we may assume that
$r:=m_G(A,z_0)>0$. First observe that it suffices to prove (\dag) only on the
circle $|\lambda|=r$. Indeed, if the inequality holds on that circle, then
by the maximum principle for subharmonic functions (applied to the function
$m_{G\times E}(A\times\{0\}, (z_0,\cdot))$) it holds for all $|\lambda|\leq r$.
In the annulus $\{r<|\lambda|<1\}$ we apply the maximum principle to the
subharmonic function $\lambda\too
\frac1{|\lambda|}m_{G\times E}(A\times\{0\},(z_0,\lambda))$.

Now fix a $\lambda_0\in E$ with $|\lambda_0|=r$. Let $f$ be an extremal function
for $m_G(A,z_0)$ with $f|_A=0$ and $f(z_0)=\lambda_0$.
Consider $F:G\too G\times E$, $F(z):=(z,f(z))$. Then
$$
m_G(A\times\{0\},(z_0,\lambda_0))\leq m_G(A,z_0)=\max\{m_G(A,z_0),\;|\lambda_0|\},
$$
which completes the proof.
\end{proof}

\noindent{\bf Acknowledgement.} The authors are indebted to professor W{\l}odzimierz Zwonek
for his valuable comments on the paper.

\halfskip

\bibliographystyle{amsplain}

\begin{thebibliography}{XXXXXXXXX}
\makeatletter\renewcommand{\@biblabel}[1]{[#1]}\makeatother
\bibitem[Car-Wie 2002]{CarWie2002}
M.~Carlehed,  J.~Wiegerinck,
\textit{Le c\^one des fonctions plurisousharmoniques n\'egatives et une
conjecture de Coman}, Ann. Polon. Math. (2002).
\bibitem[Com 2000]{Com2000a}
D.~Coman,
\textit{The pluricomplex Green function with two poles of the unit ball of $\CC^n$},
Pac. J. Math. 194 (2000), 257--283.
\bibitem[Edi 1997]{Edi1997b}
A.~Edigarian,
\textit{On the product property of the pluricomplex Green function},
Proc. Amer. Math. Soc. 125 (1997), 2855--2858.
\bibitem[Edi 1999]{Edi1999}
A.~Edigarian,
\textit{Remarks on the pluricomplex Green function},
Univ. Iag. Acta Math. 37 (1999), 159--164.
\bibitem[Edi 2001]{Edi2000b}
A.~Edigarian,
\textit{On the product property of the pluricomplex Green function, II},
Bull. Pol. Acad. Sci. 49 (2001), 389--394.
\bibitem[Edi 2002]{Edi2002}
A.~Edigarian,
\textit{Analytic discs method in complex analysis},
Diss. Math. 402 (2002), 1--56.
\bibitem[Edi-Zwo 1998]{EdiZwo1998b}
A.~Edigarian,  W.~Zwonek,
\textit{Invariance of the pluricomplex Green function under proper mappings with applications},
Complex Variables 35 (1998), 367--380.
\bibitem[Jar-Pfl 1993]{JarPfl1993}
M.~Jarnicki,  P.~Pflug,
\textit{Invariant Distances and Metrics in Complex Analysis},
 de Gruyter Expositions in Mathematics 9, Walter de Gruyter  1993.
\bibitem[L\'ar-Sig 1998]{LarSig1998b}
F.~L\'arusson,  R.~Sigurdsson,
\textit{Plurisubharmonic functions and analytic discs on manifolds},
J. Reine Angew. Math. 501 (1998), 1--39.
\bibitem[Lel 1989]{Lel1989} P.~Lelong, \textit{Fonction de Green
pluricomplexe et lemme de Schwarz dans les espaces de Banach},
J.~Math.~Pures Appl. 68 (1989), 319--347.
\bibitem[Zwo 2000a]{Zwo2000b}
W.~Zwonek,
\textit{Regularity properties of the Azukawa metric},
J.~Japan Math. Soc. 52 (2000), 899-914.
\bibitem[Zwo 2000b]{Zwo2000c}
W.~Zwonek,
\textit{Completeness, Reinhardt domains and the method of complex geode\-sics in the theory of invariant functions},
Diss. Math. 388 (2000), 1--103.

\end{thebibliography}

\end{document}